\documentclass[preprint,12pt]{elsarticle}
\usepackage{amssymb}
\usepackage{amsmath}
\usepackage{amsthm}
\usepackage{geometry}
\journal{Journal of Number Theory}
\newtheorem{defn}{Definition}
\newtheorem{thm}{Theorem}
\newtheorem{lemma}{Lemma}
\begin{document}
\begin{frontmatter}
\title{The stable property of Newton slopes for general Witt towers}

\author{Xiang Li}

\address{Morningside Center of Mathematics, Chinese Academy of Sciences, Beijing 100190. }
\ead{lixiang12@mails.ucas.ac.cn}

\begin{abstract}
Any polynomial $f(x)\in\mathbb{Z}_q[x]$ defines a Witt vector $[f]\in W(\mathbb{F}_q[x])$. Consider the Artin-Schreier-Witt tower $y^F-y=[f]$. This is a tower of curves over $\mathbb{F}_q$, with total Galois group $\mathbb{Z}_p$. We want to study the Newton slopes of zeta functions of this tower. We reduce it to the Newton polygons of L-functions associated with characters on the Galois groups. We prove that, when the conductors are large enough, these Newton slopes are unions of arithmetic progressions which are changing proportionally as the conductor increases. This is a generalization of the result of \cite{Da}, where they get the same result in the case the non-zero coefficients of $f(x)$ are roots of unity. To overcome the new difficulty in our process, we apply some $(p^{\theta},T)$-topology.
\end{abstract}


\end{frontmatter}

\section{Introduction}
Let $\mathbb{F}_q$ be a finite field of cardinality $q=p^a$ with $p$ prime. Let $\cdots\rightarrow C_n\rightarrow C_{n-1}\rightarrow\cdots\rightarrow C_0=\mathbb{P}_{\mathbb{F}_q}^1$ be a $\mathbb{Z}_p$-cover of smooth projective geometrically irreducible curves over $\mathbb{F}_q$. This means, for every $n$, $\mathrm{Gal}(C_n/C_0)\cong\mathbb{Z}/p^n\mathbb{Z}$. For the whole tower $C_{\infty}:=\lim\limits_{\longleftarrow}C_n$, we have $\mathrm{Gal}(C_\infty/C_0)\cong\mathbb{Z}_p$. Assume that each $C_n$ is totally ramified at $\infty$ and unramified outside $\infty$. By the work of M.Kosters and D.Wan \cite{Zp}, this tower corresponds to a power series $f(x)=c_0+\sum\limits_{i>0,(i,p)=1}c_ix^i\in\mathbb{Z}_q[[x]]$ which is convergent in the unit disk (i.e., $c_i\rightarrow0$ as $i\rightarrow\infty$), where $c_0\in\mathbb{Z}_p$. For every $n$, $C_n$ is the projective closure of the affine curve defined by the first $n$ equations (with respect to coordinates of Witt vectors) of
\[
y^{F}-y=c_0+\sum\limits_{i>0,(i,p)=1}c_i[x]^i,
\]
where $c_i\in\mathbb{Z}_q\cong W(\mathbb{F}_q)$, $[x]=(x,0,0,\ldots),y=(y_1,y_2,\ldots)$ are $p$-typical Witt vectors, and $\cdot^F$ is the Frobenius map raising every coordinate of a Witt vector to its $p$-th power. The additions and multiplications in both sides are those of Witt vectors. Without loss of generality, one can assume $c_0=0$.

Consider the zeta function of $C_n$. By Weil's theorems \cite{Weil}, \[\zeta(C_n,s)=\exp\left(\sum\limits_{k\geq1}\# C_n(\mathbb{F}_{q^k})\cdot\frac{s^k}{k}\right)=\frac{P(C_n,s)}{(1-s)(1-qs)}.\] $P(C_n,s)$ is a polynomial of degree $2g_n$, where $g_n$ is the genus of $C_n$. In the spirit of Iwasawa theory, it is natural to study the $q$-adic Newton polygons of the sequence $P(C_n,s)$, especially the behavior of them when $n\rightarrow\infty$.

$P(C_n,s)$ is a product of $L$-functions for various finite characters. Let $\chi:\mathbb{Z}_p\rightarrow\mathbb{C}_p^{\times}$ be a finite additive character. That means the image of $\chi$ has finite cardinality. We call this cardinality the conductor of $\chi$. Let $m_{\chi}$ be the nonnegative integer such that the conductor of $\chi$ is $p^{m_{\chi}}$. If $\chi$ is nontrivial, we let $\pi_{\chi}=\chi(1)-1$, then $\nu_p(\pi_{\chi})=\varphi(p^{m_{\chi}})^{-1}$, where $\varphi$ is the Euler function. Note that $\mathrm{Gal}(C_{m_{\chi}}/C_0)\cong\mathbb{Z}/p^{m_{\chi}}\mathbb{Z}$, so $\chi$ can be seen as a homomorphism $\mathrm{Gal}(C_{m_{\chi}}/C_0)\rightarrow\mathbb{C}_p^{\times}$. Now, we define the $L$-function $L(\chi,s)$ over $\mathbb{A}_{\mathbb{F}_q}^1$ by
\[L(\chi,s)=\prod_{x\in|\mathbb{A}_{\mathbb{F}_q}^1|}\frac{1}{1-\chi(\mathrm{Frob}_x)s^{\deg(x)}},\] where $|\mathbb{A}^1_{\mathbb{F}_q}|$ is the set of closed points of $\mathbb{A}^1_{\mathbb{F}_q}$, and $\mathrm{Frob}_x\in\mathrm{Gal}(C_{m_{\chi}}/C_0)$ is the Frobenius element associated with $x$. By Weil's result \cite{Weil}, $L(\chi,s)$ is also a polynomial. In fact, $L(\chi,s)\in1+s\mathbb{Z}_p[\pi_{\chi}][s]$.

For different choices of $\chi$ with fixed $m_{\chi}$, the polynomials $L(\chi,s)$ are conjugated with each other. Denote the degree of $L(\chi,s)$ by $d_{m_{\chi}}$. After decomposing in $\mathbb{C}_p[s]$ and adjusting the order of factors, we have $L(\chi,s)=\prod_{i=1}^{d_{m_{\chi}}}(1-\alpha_i^{(\chi)}s)$ such that $0<\nu_q(\alpha_1^{(\chi)})\leq\ldots\leq\nu_q(\alpha_{d_{m_{\chi}}}^{(\chi)})<1$. For fixed $m_{\chi}$, the $q$-slope sequence $\{\gamma_i^{(m_{\chi})}:=\nu_q(\alpha_i^{(\chi)})\}_{1\leq i\leq d_{m_{\chi}}}$ depends only on $m_\chi$, not on $\chi$.

It is a standard fact that \[P(C_n,s)=\prod\limits_{1\leq m\leq n}\prod\limits_{\mathrm{cond}(\chi)=p^m}L(\chi,s).\] Thus, the study of the polynomial $P(C_n,s)$ reduces to the study of $L(\chi,s)$ for various finite characters $\chi$, or equivalently, the slope sequences $\{\gamma_i^{(m)}\}_{1\leq i\leq d_m}$ for various $m$.

\noindent\textbf{Notation.} For a discrete valuation ring $A$ with nontrivial valuation $\nu$, the Newton polygon of $\sum_{n\geq0} a_nx^n\in A[[x]]$ with respect to $\nu$ is the lower convex hull of the set of points $\{(n,\nu(a_n))|n\geq0\}$.

\noindent\textbf{Unit Root Case.} Denote the set of Teichm\"uller liftings of $\mathbb{F}_q$ by $\widehat{\mathbb{F}}_q$, and the set of polynomials with coefficients in $\widehat{\mathbb{F}}_q$ by $\widehat{\mathbb{F}}_q[x]$. In 2013, Christopher Davis, Daqing Wan and Liang Xiao \cite{Da} proved that for $f(x)\in\widehat{\mathbb{F}}_q[x]$, the corresponding tower $C_{\infty}$ is slope stable with respect to the following sense.

\begin{defn} The $\mathbb{Z}_p$-tower $C_{\infty}$ is called slope stable, if there exists a positive integer $m_0$ such that for $m\geq m_0$, the $q$-slope sequence
\[
\{\gamma_i^{(m)},1\leq i\leq d_m\}=\bigcup\limits_{j=0}^{p^{m-m_0}-1}\left\{\frac{j}{p^{m-m_0}},\frac{\gamma_1^{(m_0)}+j}{p^{m-m_0}},\ldots,\frac{\gamma_{d_{m_0}}^{(m_0)}+j}{p^{m-m_0}}\right\}\backslash\{0\}.
\]
\end{defn}

\smallskip

Now, we state our main theorem, which is a generalization of the above result.
\begin{thm}[Polynomial Case]\label{main}
For $f(x)=c_0+\sum\limits_{(i,p)=1, D\geq i>0}c_ix^i\in\mathbb{Z}_q[x]$, the corresponding tower $C_{\infty}$ is slope stable.
\end{thm}

In the rest of this paper, we give a proof of this theorem, finding the appropriate $m_0$. The method is analogous to \cite{Da}. First, we consider the equivalent characteristic functions built out of $L$-functions. Then, we construct the $T$-adic characteristic function, which interpolates every characteristic function of nontrivial finite conductor. We prove the Dwork's trace formula of $T$-adic character functions; the calculation process is analogous to \cite{Wan}. After that, we produce a lower bound of the Newton polygon of the $T$-adic characteristic function and an upper bound of the $\pi_{\chi}$-adic Newton polygon of the characteristic function for a finite character $\chi$ with a specific conductor. The two bounds touch periodically. Then, we use the isolate properties of valuations of $\mathbb{Z}_p$ to prove that the Newton polygons for finite characters have a stable property if the conductor is large enough.

The new difficulty, comparing with \cite{Da}, is, in our case, if we apply $T$-adic valuation to the $T$-adic characteristic function, the $p$-power factor in the coefficients will be ignored. Then, we fail to produce a sharp enough lower bound of the $\pi_{\chi}$-adic Newton polygon to touch the upper bound. The solution is, we consider the finer $(p^{\theta},T)$-adic valuation for a suitable $\theta$ instead of $T$-adic valuation, and the computation of the inequality is a bit more tricky. By these modifications, we take the $p$-power factor into acount, and get over the difficulty.

Recently, R.Ren, D.Wan, L.Xiao and M.Yu \cite{Nw} proved a slope stable property in unit root case for higher rank (i.e., the $\mathbb{Z}_p^l$-towers for $l>1$), which is a different generalization of the main result of \cite{Da}. So, it will be natural to give a conjecture that there is a slope stable property for higher rank in our polynomial case.

\smallskip

\noindent\textbf{Acknowledgements.} I would like to thank Professor Daqing Wan for his critical advice when he was visiting the Morningside Center of Mathematics. I also thank Professor Qifang Zhang and my supervisor Professor Ye Tian for introducing this topic to me and for their constant help.
\section{Slope stability with respect to characteristic functions}
To simplify the definition of slope stability, we introduce the equivalent characteristic functions which are built out of $L$-functions.

First, we drop the point $0\in|\mathbb{A}^1_{\mathbb{F}_q}|$, and consider the $L$-function defined on the torus. For a finite character $\chi$, let
\[L^*(\chi,s)=\prod\limits_{x\in|\mathbb{G}_m|}\frac{1}{1-\chi(\mathrm{Frob}_x)s^{deg(x)}}=(1-\chi(\mathrm{Frob}_0)s)L(\chi,s).\]
Then, we define the characteristic function \[C^*(\chi,s)=\prod\limits_{j=0}^{\infty}L^*(\chi,q^js).\] It is a trivial fact that this infinite product is convergent.

Now we simplify the definition of slope stability, using $C^*$. We have the following lemma.
\begin{lemma}\label{conv}
The $\mathbb{Z}_p$-tower $C_{\infty}$ is slope stable with respect to a positive integer $m_0$, if the $\pi_{\chi}$-adic Newton polygon of the corresponding $C^*(\chi,s)$ is the same for any $\chi$ such that $m_{\chi}\geq m_0$, where $\pi_{\chi}=\chi(1)-1$.
\end{lemma}
\begin{proof}
Let $\chi_0$ be a character such that $m_{\chi_0}=m_0>0$. The $q$-adic Newton slopes of $L(\chi_0,s)$ are $\gamma_1^{(m_0)}<...<\gamma_{d_{m_0}}^{(m_0)}$, satisfying $0<\gamma_j^{(m_0)}<1$, $\forall j$. As $L^*(\chi_0,s)=(1-\chi_0(\mathrm{Frob}_0)s)L(\chi_0,s)$ and $\chi_0$ maps any element in $\mathbb{Z}_p$ into a root of unity, the $q$-adic Newton slopes of $L^*(\chi,s)$ are $0,\gamma_1^{(m_0)},...,\gamma_{d_{m_0}}^{(m_0)}$. By definition, \[C^*(\chi_0,s)=\prod\limits_{j=0}^{\infty}L^*(\chi_0,q^js),\] hence the set of $q$-adic Newton slopes of $C^*(\chi_0,s)$ is \[\bigcup\limits_{i\geq0}\left\{i,\gamma_1^{(m_0)}+i,...,\gamma_{d_{m_0}}^{(m_0)}+i\right\}.\] Since $\nu_q(\pi_{\chi_0})=\frac{1}{a\varphi(p^{m_0})}$, the set of $\pi_{\chi_0}$-adic Newton slopes of $C^*(\chi_0,s)$ is \[\bigcup\limits_{i\geq0}\{ia\varphi(p^{m_0}),(\gamma_1^{(m_0)}+i)a\varphi(p^{m_0}),...,(\gamma_{d_{m_0}}^{(m_0)}+i)a\varphi(p^{m_0})\}.\] Here $\varphi$ is the Euler function. By the assumption of this lemma, these values are also the $\pi_{\chi}$-adic Newton slopes of $C^*(\chi,s)$ for any $\chi$ such that $m_{\chi}\geq m_0$. As $\nu_q(\pi_{\chi})=\frac{1}{a\varphi(p^{m_{\chi}})}$, and $\frac{a\varphi(p^{m_0})}{a\varphi(p^{m_{\chi}})}=\frac{1}{p^{m_{\chi}-m_0}}$, the $q$-adic Newton slopes of $C^*(\chi,s)$ are \[\bigcup\limits_{i\geq0}\left\{\frac{i}{p^{m_{\chi}-m_0}},\frac{\gamma_1^{(m_0)}+i}{p^{m_{\chi}-m_0}},...,\frac{\gamma_{d_{m_0}}^{(m_0)}+i}{p^{m_{\chi}-m_0}}\right\}.\] In view of the relation
\[L(\chi,s)=\frac{1}{1-\chi(\mathrm{Frob}_0)s}\frac{C^*(\chi,s)}{C^*(\chi,qs)},\] the $q$-adic Newton polygon of $L(\chi,s)$
has slopes \[\bigcup\limits_{i=0}^{p^{m_{\chi}-m_0}-1}\left\{\frac{i}{p^{m_{\chi}-m_0}},
\frac{\gamma_1^{(m_0)}+i}{p^{m_{\chi}-m_0}},...,\frac{\gamma_{d_{m_0}}^{(m_0)}+i}{p^{m_{\chi}-m_0}}\right\}-\{0\}.\qedhere\]
\end{proof}

Therefore, the main problem of this paper comes down to the problem of the stability of $\pi_{\chi}$-adic Newton polygons of $C^*(\chi,s)$ for various $\chi$ with large enough conductor.

\section{Exponential sums and $T$-adic functions}
In this section, we will do some preparations. First, we introduce the exponential form of the $L$-functions, which is much more convenient for calculation. The deduction is due to the work of Chunlei Liu and Dasheng Wei \cite{Wei}.

Assume that $f(x)\in\mathbb{Z}_q[x]$ is a polynomial of degree $D$ whose nonconstant terms $c_ix^i$ satisfy $(i,p)=1$. Write $f(x)=\sum\limits_{j=0}^{\infty}p^jf_j(x)$, where $f_j(x)\in\widehat{\mathbb{F}}_q[x]$ has degree $\widetilde{d_j}$, $\forall j$. Let $m$ be a positive integer, and $\chi$ be a character such that $m_{\chi}=m$. Let $f^{(m)}(x)=\sum\limits_{j=0}^{m-1}p^jf_j(x)$. By Lemma 1.1 of \cite{Wei}, we have \[L^*(\chi, s) = \exp(\sum\limits_{k=1}^{\infty}S^*(\chi, k)\frac{s^k}{k}),\] where
\[S^*(\chi, k)=\sum\limits_{x\in\mu_{q^k-1}}\chi(\mathrm{Tr}_{\mathbb{Q}_{q^k}/\mathbb{Q}_p}(f^{(m)}(x))).\]
Note that $\mathrm{Cond}(\chi)=p^m$, so $\chi(\mathrm{Tr}_{\mathbb{Q}_{q^k}/\mathbb{Q}_p}(p^m\lambda))=1$ for any $\lambda\in\mathbb{Z}_{q^k}$. Therefore, \[S^*(\chi, k)=\sum\limits_{x\in\mu_{q^k-1}}\chi(\mathrm{Tr}_{\mathbb{Q}_{q^k}/\mathbb{Q}_p}(f(x))).\]

Next, in an analogous manner with \cite{Wan}, we define the $T$-adic $L$-function. This function can be seen as a "universal" $L$-function.

For a positive integer $k$, we define the $T$-adic exponential sum of $f$ over $\mathbb{F}_{q^k}$ as the sum
\[S^*(T, k) = \sum_{x \in \mu_{q^k-1}}(1+T)^{\mathrm{Tr}_{\mathbb{Q}_{q^k}/\mathbb{Q}_p}(f(x))}.\]
Then, the $T$-adic L-function is defined as \[L^*(T, s) = \exp(\sum\limits_{k=1}^{\infty}S^*(T, k)\frac{s^k}{k}) \in 1 + s\mathbb{Z}_p[[T]][[s]].\]
We also define the $T$-adic characteristic function $C^*(T,s)$ as
\[C^*(T,s)=\prod\limits_{j=0}^{\infty}L^*(T,q^js)=\exp(-\sum\limits_{k=1}^{\infty}\frac{1}{q^k-1}S^*(T,k)\frac{s^k}{k}).\]
Conceptually, the $T$-adic $L$-function is defined in a similar manner with the original $L$-functions involving finite characters. In fact, we just need to replace the finite character $\chi$ by
\begin{equation*}
\chi_T:\mathbb{Z}_p\rightarrow\mathbb{Z}_p[[T]]^{\times},1\mapsto 1+T.
\end{equation*}
The relations between them is clear: $L^*(T, s)|_{T = \pi_{\chi}} = L^*(\chi, s),$ and $C^*(T, s)|_{T = \pi_{\chi}} = C^*(\chi, s).$

Finally, we take a variable substitution, which is a preparation for further computations in the next section. Consider the Artin-Hasse exponential series
$E(X) = \exp(\sum\limits_{i=1}^{\infty}\frac{X^{p^i}}{p^i})\in 1+X+X^2\mathbb{Z}_p[[X]].$ Let $t=E(X)-1$, then simple iteration implies $X \in t+t^2\mathbb{Z}_p[[t]]$. We introduce new indeterminates $\pi$ and $\pi_i$ ($i\geq0$). Let $\pi$, $\pi_i$ be power series in $\mathbb{Z}_p[[T]]$ such that $E(\pi) = 1+T$, and $E(\pi_i) = (1+T)^{p^i}$ for any $i \geq 0$. We have $\pi_0 = \pi$. By the properties of Artin-Hasse exponential
series given above, we have $T \in \pi+\pi^2\mathbb{Z}_p[[\pi]]$ and $\pi \in T+T^2\mathbb{Z}_p[[T]]$. As $E(\pi_{i+1})=E(\pi_i)^p$, we have $\pi_{i+1}\in p\pi_i\mathbb{Z}_p[[\pi]]+\pi_i^p\mathbb{Z}_p[[\pi]]$. By induction, $\pi_i\in\sum\limits_{j=0}^{i} p^{i-j}\pi^{p^j}\mathbb{Z}_p[[\pi]], \forall i$.

We construct the function $E_f$, which plays, roughly speaking, the role of Dwork's $G$-function in our story. For any $h(x)=\sum\limits_{u=0}^{d'}h_ux^u \in \widehat{\mathbb{F}}_q[x]$, and any indeterminate $\eta$ independent with $x$, define
\[E_{h}(x,\eta)=\prod\limits_{u=0}^{d'}E(\eta{h_u}x^u).\] Let $\sigma\in\mathrm{Gal}(\mathbb{Q}_q/\mathbb{Q}_p)$ be the Frobenius endomorphism such that $\sigma(x)=x^p$ for all
$x\in\mu_{q-1}$, and let it act on $\mathbb{Q}_q[[\pi]]$, leaving $\pi$ fixed. Remember that $f(x)=\sum\limits_{j=0}^{\infty}p^jf_j(x)$, where $f_j(x)\in\widehat{\mathbb{F}}_q[x]$, for every $j$. Define $E_f(x)=\prod\limits_{j=0}^{\infty}E_{f_j}(x,\pi_j)$. The convergence of this infinite product is an easy consequence of the fact that $\pi_i\in\sum\limits_{j=0}^{i} p^{i-j}\pi^{p^j}\mathbb{Z}_p[[\pi]]$. We have the Dwork's splitting formula (Lemma 4.3 of \cite{Wan}):
\begin{lemma}
If $x \in \mu_{q-1}$, then
\[E(\eta)^{\mathrm{Tr}_{\mathbb{Q}_{q}/\mathbb{Q}_p}(h(x))}=\prod\limits_{i=0}^{a-1}E_h^{\sigma^i}(x^{p^i},\eta).\]
\end{lemma}

\begin{proof}
The proof is given in section 4 of \cite{Wan}.
\end{proof}

Now, for $x\in\mu_{q^k-1}$,
\begin{eqnarray*}
(1+T)^{\mathrm{Tr}_{\mathbb{Q}_{q^k}/\mathbb{Q}_p}(f(x))}&=&\prod\limits_{j=0}^{\infty}(1+T)^{p^j\mathrm{Tr}_{\mathbb{Q}_{q^k}/\mathbb{Q}_p}(f_j(x))}
=\prod\limits_{j=0}^{\infty}E(\pi_j)^{\mathrm{Tr}_{\mathbb{Q}_{q^k}/\mathbb{Q}_p}(f_j(x))}\\&=&\prod\limits_{j=0}^{\infty}\prod_{i=0}^{ak-1}E_{f_j}^{\sigma^i}(x^{p^i},\pi_j)
=\prod\limits_{i=0}^{ak-1}E_{f}^{\sigma^i}(x^{p^i}).
\end{eqnarray*}

As a conclusion, we have

\begin{lemma}
If $x \in \mu_{q^k-1}$, then \[(1+T)^{\mathrm{Tr}_{\mathbb{Q}_{q^k}/\mathbb{Q}_p}(f(x))}=\prod\limits_{i=0}^{ak-1}E_f^{\sigma^i}(x^{p^i}).\]
\end{lemma}

Therefore, the $T$-adic exponential sum $S^*(T,s)$ has the following new expression:
\[S^*(T,s)=\sum_{x \in \mu_{q^k-1}}\prod\limits_{i=0}^{ak-1}E_f^{\sigma^i}(x^{p^i}).\]
\section{Dwork's trace formula}
We interpret $C^*(T,s)$ as the characteristic function of a compact operator, proving its analytic continuation. It is accustomed to call this type of result ``Dwork's trace formula'', since the basic method was developed by Bernard Dwork. In a similar method, we extend Chunlei Liu and Daqing Wan's work \cite{Wan} from the unit root case to more general polynomial case, where we meet a new difficulty to deal with the $p$-orders of the coefficients. In our case, the contribution of $p$-orders must be counted, so we have to replace $T$-adic valuation (which is used in \cite{Wan}) by $(p^{\theta},T)$-adic valuation for a selected $\theta>0$.

Recall that $f(x)=\sum\limits_{j\geq0}p^jf_j(x)$ is a polynomial of degree $D$, where $\mathrm{deg}(f_j)=\widetilde{d_j}$ for every $j$. Let $d=\max\{\frac{\widetilde{d_j}}{p^j}|j\geq0\}$. Let $l$ be the first $j$ such that $D=\widetilde{d_j}$.

We will define a ring $\Lambda$ with a valuation on it, which will be the coefficient ring for the space where the compact operator acts.

First, assume $l>0$. Choose a positive rational number $\theta$ such that $\theta<\frac{1}{\varphi(p^{l})}$. Let $N$ be a positive integer such that $Nd^{-1}$, $N{\theta}^{-1}$ are integers. Define a ring \[\Lambda=\mathbb{Z}_q[[\pi^{N^{-1}},\pi',\pi'^{-1}]]/(p\pi'-\pi^{\theta^{-1}})\] and a $\Lambda$-module $R=(p^{\theta},\pi)$. Now, we define the $R$-adic valuation $\nu_R$. Let $\nu_R(\lambda \pi^j)=\theta^{-1}\nu_p(\lambda)+j$ for any $\lambda\in\mathbb{Z}_q$ and $j\geq0$. If $\omega$ is the valuation defined by the ideal $(p,\pi^{\theta^{-1}})$ of $\Lambda$, then $\nu_R=\theta^{-1}\omega$.

If $l=0$, set $\theta=+\infty$, $\Lambda=\mathbb{Z}_q[[\pi]]$, $R=(\pi)$ and $\nu_R=\nu_{\pi}$. This is the degenerate case.

We find a lower bound of $\nu_R(\pi_j)$ for any $j\geq0$.

\begin{lemma}\label{Rvalue}
Let $l$, $\theta$, $\Lambda$ and $R$ be defined as above, then we have $\nu_R(\pi_j)\geq p^j$ for $0\leq j\leq l$, and $\nu_R(\pi_j)\geq p^{l}$ for $j\geq l$.
\end{lemma}

\begin{proof}
In the $\theta=\infty$ case, this lemma is just saying $\nu_{\pi}(\pi)\geq1$. Now we assume that $\theta<\infty$. For $j=0$, the inequality is automatically true. When $j<l$, if $\nu_R(\pi_j)\geq p^j$, then $\nu_R(\pi_j^p)\geq p^{j+1}$, and $\nu_R(p\pi_j)\geq \frac{1}{\theta}+p^j\geq\varphi(p^l)+p^j\geq\varphi(p^{j+1})+p^j=p^{j+1}$. Hence, the first part of the lemma is proved. When $j\geq l$, we have $\nu_R(\pi_{j+1})\geq\nu_R(\pi_j)$. This implies the second part.
\end{proof}

Write $E_f(x)=\sum\limits_{u\geq0}\alpha_ux^u$, where $\alpha_u\in\mathbb{Z}_q[[\pi]], \forall u$. We study the $R$-adic valuation of $\alpha_u$.

\begin{thm}\label{lower}
$\nu_R(\alpha_u)\geq\frac{u}{d}$, $\forall u\geq0$.
\end{thm}
\begin{proof}
Review the definition: $E_f(x)=\prod\limits_{j=0}^{\infty}E_{f_j}(x,\pi_j)$. Define $\alpha_u^{(j)}$ such that $E_{f_j}(x,\pi_j)=\sum\limits_{u\geq0}\alpha_u^{(j)}x^u$. Write $f_j(x)=\sum\limits_{i=0}^{\widetilde{d_j}}c_{ij}x^i$. By definition, $E_{f_j}(x,\pi_j)=\prod\limits_{i}E(\pi_jc_{ij}x^i)$. As a consequence, we have $\nu_{\pi_j}(\alpha_u^{(j)})\geq\frac{u}{\widetilde{d_j}}$ as long as $\widetilde{d_j}>0$. Thus, $\nu_R(\alpha_u^{(j)})\geq\frac{u}{\widetilde{d_j}}\cdot\nu_R(\pi_j)$. From Lemma \ref{Rvalue}, we have $\nu_R(\alpha_u^{(j)})\geq\frac{u}{\frac{\widetilde{d_j}}{p^j}}$ for $0\leq j\leq l$, and $\nu_R(\alpha_u^{(j)})\geq\frac{u}{\frac{\widetilde{d_j}}{p^l}}$ for $j\geq l$. Note that $E_f(x)=\sum\limits_{u\geq0}\alpha_ux^u=\prod\limits_{j=0}^{\infty}\sum\limits_{u\geq0}\alpha_u^{(j)}x^u$. As a conclusion, we have \[\nu_R(\alpha_u)\geq\frac{u}{\max\limits_{j\geq0}\{\{\frac{\widetilde{d_j}}{p^j}|0\leq j\leq l\}\bigcup\{\frac{\widetilde{d_j}}{p^l}|j\geq l\}\}}=\frac{u}{d}.\qedhere\]
\end{proof}

Next, we define the appropriate space and operator for the trace formula of $C^*(T,s)$. Let $B$ be the Banach $\Lambda$-module with orthonormal basis $\{Y_u=\pi^{\frac{u}{d}}x^u|u\geq0\}$, and $B'$ be the module with the same set as formal basis. To be precise, \[B'=\{\sum\limits_{u\geq0}\beta_uY_u|\beta_u\in\Lambda\},\]\[B=\{\sum\limits_{u\geq0}\beta_uY_u\in B'|\nu_R(\beta_u)\rightarrow\infty\ as\ u\rightarrow\infty\}.\] In the above setting, Theorem \ref{lower} is equivalent to the fact that $E_f(x)\in B'$.

Let $\psi_p$ be the map $\sum\limits_{u\geq0}\beta_uY_u\mapsto\sum\limits_{u\geq0}\beta_{pu}\pi^{\frac{(p-1)u}{d}}Y_u$ from $B'$ to $B$. Let $\sigma\in\mathrm{Gal}(\mathbb{Q}_q/\mathbb{Q}_p)$ be the Frobenius map acting on $\mathbb{Q}_q$, fixing $\pi$ and $x$. Define $\psi=\sigma^{-1}\circ\psi_p\circ E_f(x)$, where $E_f(x)$ means the map multiplying by $E_f(x)$ from $B$ to $B'$. Thus, $\psi$ maps from $B$ to $B$. Note that $\psi_p\circ G(x^p)=G(x)\circ\psi_p$ for any $G(x)\in B'$, we have \begin{eqnarray*}\psi^{r}(g)&=&\psi^{r-1}(\psi_p(E_f^{\sigma^{-1}}(x)g^{\sigma^{-1}}))
=\psi^{r-2}(\psi_p^2(E_f^{\sigma^{-1}}(x^p)E_f^{\sigma^{-2}}(x)g^{\sigma^{-2}}))\\&=&...
=\psi_p^r(g^{\sigma^{-r}}\prod\limits_{i=0}^{r-1}E_f^{\sigma^{i-r}}(x^{p^i}))
\end{eqnarray*}
for $r\geq0$, $g\in B$. Let $r=ak$, we find that
\[\psi^{ak}=\psi_p^{ak}\circ \prod\limits_{i=0}^{ak-1}E_f^{\sigma^i}(x^{p^i}),\]
where $\prod\limits_{i=0}^{ak-1}E_f^{\sigma^i}(x^{p^i})$ means multiplying by $\prod\limits_{i=0}^{ak-1}E_f^{\sigma^i}(x^{p^i})$. $\psi^{ak}$ is linear over $\Lambda$ for any $k\geq1$, but $\psi$ is only semi-linear over $\Lambda$.

Write $\prod\limits_{i=0}^{ak-1}E_f^{\sigma^i}(x^{p^i})=\sum\limits_{u\geq0}\beta_{u,k}x^u$, where $\beta_{u,k}\in\Lambda$. Define $\beta_{u,k}=0$ for any $u<0$.

Using the formulas above, $\psi^{ak}$ sends $Y_v$ to $\sum\limits_{u}\beta_{q^ku-v,k}\pi^{\frac{v-u}{d}}Y_u$. For $k=1$, $M=(\beta_{qu-v,1}\pi^{\frac{v-u}{d}})_{u,v}$ is the matrix of $\psi^a$ over the basis $\{Y_u|u\geq0\}$. By Theorem \ref{lower} and the equation $\prod\limits_{i=0}^{ak-1}E_f^{\sigma^i}(x^{p^i})=\sum\limits_{u\geq0}\beta_{u,k}x^u$, we have \[\nu_R(\beta_{u,k})\geq\frac{u}{p^{ak-1}d}.\] We calculate that $\frac{qu-v}{p^{a-1}d}+\frac{v-u}{d}=\frac{(p-1)u}{d}+\frac{(q-p)v}{qd}$. Thus, the $R$-valuations of the elements of the matrix $M=(\beta_{qu-v,1}\pi^{\frac{v-u}{d}})_{u,v}$ increase at least linearly as the row index $u\rightarrow\infty$. Therefore, $M$ is nuclear. Equivalently, $\psi^a$ is a compact operator on the Banach $\Lambda$-module $B$. Thus, for any $k\geq1$, we can consider the nuclear trace of $\psi^{ak}$. In fact, it is just $\mathrm{Tr}(M^k)$.

Now, we are ready to prove the Dwork's trace formula in our case.

\begin{thm}[Dwork's trace formula]\label{dwork}
For every positive integer $k$,
\[S^*(T,k)=\sum\limits_{x\in\mu_{q^k-1}}\prod\limits_{i=0}^{ak-1}E_f^{\sigma^i}(x^{p^i})=(q^k-1)\mathrm{Tr}(\psi^k|B/\Lambda).\]
\end{thm}
\begin{proof}
\begin{align*}
S^*(T,k)&=\sum\limits_{x\in\mu_{q^k-1}}\prod\limits_{i=0}^{ak-1}E_f^{\sigma^i}(x^{p^i})
=\sum\limits_{x\in\mu_{q^k-1}}\sum\limits_{u\geq0}\beta_{u,k}x^u\\
&=\sum\limits_{u\geq0}(\sum\limits_{x\in\mu_{q^k-1}}x^u)\beta_{u,k}=(q^k-1)\sum\limits_{u\geq0}\beta_{(q^k-1)u,k}\\
&=(q^k-1)\mathrm{Tr}(M^k)=(q^k-1)\mathrm{Tr}(\psi^{ak}|B/\Lambda).\qedhere
\end{align*}
\end{proof}

Now, we have the following description of $C^*(T,s)$, which is the multiplicative form of the Dwork's trace formula.

\noindent\textbf{Notation.} Generally, we say a power series $\sum\limits_{n\geq0}a_nx^n$ is analytic (with respect to some valuation), if it converges for all $x$ in the definition domain.
\begin{thm}[Analytic trace formula]\label{dwork2}
$$C^*(T,s)=\det(1-s\psi^a|B/\Lambda),$$ and the function $C^*(T,s)$ is $R$-adic analytic.
\end{thm}
\begin{proof}
$\det(1-s\psi^a|B/\Lambda)$ is just $\det(I-sM)$. This formula is a special case of Lemma 4.2 of \cite{Wan1}. The existence of the determinant $\det(I-sM)$ is derived from the facts that $\nu_R(\beta_{u,1})\geq\frac{u}{p^{a-1}d}$ and $M$ has nuclear shape such that the $R$-adic valuation of the elements of $M$ increase at least linearly as the row index tends to $\infty$. Using Dwork's trace formula (theorem \ref{dwork}) and the well-known formula (for a proof, see \cite{Wan1})\[\det(I-sA)=\exp(-\sum\limits_{k=1}^{\infty}\mathrm{Tr}(A^k)\frac{s^k}{k}),\] we have
\begin{eqnarray*}
C^*(T,s)&=&\exp(-\sum\limits_{k=1}^{\infty}\frac{1}{q^k-1}S^*(T,k)\frac{s^k}{k})\\
&=&\exp(-\sum\limits_{k=1}^{\infty}\mathrm{Tr}(M^k)\frac{s^k}{k})=\det(I-sM).
\end{eqnarray*}
By the nuclear shape of $M$ again, $C^*(T,s)$ is $R$-adic analytic.
\end{proof}
\section{Lower bound of the $R$-adic Newton polygon}
Now, from Theorem \ref{lower} and Theorem \ref{dwork2}, we find a lower bound of the $R$-adic Newton polygon of $C^*(T,s)$.

\begin{thm}\label{lower2}
The $R$-adic Newton polygon of $C^*(T,s)$ lies above (maybe coincides with) the polygonal line connecting $\{(n,\frac{a(p-1)n(n-1)}{2d})|n\geq0\}$.
\end{thm}
\begin{proof}
Note that $\mathbb{Q}_q$ can be obtained by adding a $(p^a-1)$-th root of unity to $\mathbb{Q}_p$, so there is a normal integral basis of $\mathbb{Q}_q$ over $\mathbb{Q}_p$, say $\{\xi_0,...,\xi_{a-1}\}$. Let $\Lambda_1=\mathbb{Z}_p[[\pi^{N^{-1}},\pi',\pi'^{-1}]]/(p\pi'-\pi^{\theta^{-1}})$, then $\{\xi_iY_u|0\leq i\leq{a-1},u\geq0\}$ is a set of normal basis of $B$ over $\Lambda_1$. $\psi$ is linear over $\Lambda_1$. Arrange this basis in the lexicographical order $\{\xi_0,...,\xi_{a-1},\xi_0Y_1,...\xi_{a-1}Y_1,...,\xi_0Y_u,...,\xi_{a-1}Y_u,...\}$. We compute the matrix of $\psi$ over this basis. We have \[\psi(\xi_jY_v)=\sum\limits_{u=0}^{\infty}\sigma^{-1}(\xi_j\alpha_{pu-v})\pi^{\frac{v-u}{d}}Y_u.\] Let \[\sigma^{-1}(\xi_j\alpha_{pu-v})=\sum\limits_{i=0}^{a-1}\alpha_{(i,u),(j,v)}\xi_i,\] where $\alpha_{(i,u),(j,v)}\in\Lambda_1$, we have
\[\psi(\xi_jY_v)=\sum\limits_{u=0}^{\infty}\sum\limits_{i=0}^{a-1}\alpha_{(i,u),(j,v)}\pi^{\frac{v-u}{d}}\xi_iY_u.\] Hence, the matrix of $\psi$ under the above basis is $M_1=(\alpha_{(i,u),(j,v)}\pi^{\frac{v-u}{d}})_{(i,u),(j,v)}$. By definition, if $pu<v$, then $\alpha_{(i,u),(j,v)}=0$. Therefore, powers of $M_1$ make sense, and the matrix of $\psi^r$ is $M_1^r$ for any $r\geq1$.

Consider the norm of $\Lambda((s))$ over $\Lambda_1((s))$. By simple linear algebra (the same argument with the proof of \cite[Theorem 4.8]{Wan}), we have
\[\mathrm{Norm}(\det(1-s^a\psi|B/\Lambda))=\det(1-s^a\psi|B/\Lambda_1)=\prod_{\zeta^a=1}\det(1-\zeta{s}\psi|B/\Lambda_1),\] in other words,
\[\mathrm{Norm}(\det(I-s^aM))=\det(I-s^aM_1^a)=\prod_{\zeta^a=1}\det(I-\zeta{s}M_1),\] so the $R$-adic Newton polygons of $\det(I-s^aM)$ and $\det(I-sM_1)$ coincide.

Note that for any $\omega\in\mathbb{Z}_q$, $\omega=\sum\limits_{i=0}^{a-1}\omega_i\xi_i$, if $\nu_{p}(\omega)\geq{\lambda}$, the same will be
true for every $\omega_i$. By theorem \ref{lower} and the equation \[\sigma^{-1}(\xi_j\alpha_{pu-v})=\sum\limits_{i=0}^{a-1}\alpha_{(i,u),(j,v)}\xi_i,\] we have \[\nu_R(\alpha_{(i,u),(j,v)})\geq\frac{pu-v}{d}.\] By the shape of $M_1$, the $R$-order of coefficient of $s^{an}$ of $\det(I-sM_1)$ is at least
\[0+\frac{a(p-1)}{d}+\frac{2a(p-1)}{d}+\ldots+\frac{(n-1)a(p-1)}{d}=\frac{a(p-1)n(n-1)}{2d}.\] Thus, the polygonal line connecting $$\{(an,\frac{a(p-1)n(n-1)}{2d})|n\geq0\}$$ is a lower bound of the $R$-adic Newton polygon of $\det(I-s^aM)$. But we have $C^*(T,s)=\det(1-sM)$ (Theorem \ref{dwork2}), so the polygonal line connecting $$\{(n,\frac{a(p-1)n(n-1)}{2d})|n\geq0\}$$ is a lower bound of the $R$-adic Newton polygon of $\det(I-sM)$.
\end{proof}
\section{Periodicity of Newton polygons}
In this section, we prove the main theorem. By Lemma \ref{conv}, we only need to study the Newton polygon of $C^*(\chi,s)$. In the proof of the main theorem in \cite{Da}, there is an important trick transferring between Newton polygons of $C^*(\chi,s)$ and $C^*(T,s)$. We explain this trick by the following lemma. We will do some modifications to fit our case.
\begin{lemma}\label{trick}

Consider $D(T,s)\in\mathbb{Z}_p[[T]][[s]]$. Let $H,H'$ be two infinite convex polygons starting with $(0,0)$ satisfying
\begin{itemize}
\item $H'\geq H$ with $\Delta:=\sup_{t\geq 0} (H'(t)-H(t))<\infty$;
\item there are integers $h>h'>0$ such that for all $k\geq 0$,  $H=H'$ on $[kh, kh+h']$,  $H'$ is a line on $[kh+h', (k+1)h]$,
and both $kh$ and $kh+h'$ are vertices of $H$.
\end{itemize}
Assume that there is an integer $m\geq1$ and a character
$\chi_1$ with $m_{\chi_1}=m$. We also assume that there is a positive rational number $\theta>\varphi(p^m)^{-1}$ such that $D(T,s)$ is $(p^{\theta},T)$-adic entire, and $$H\leq \mathrm{NP}_{(p^\theta, T)}(D(T,s))\leq \mathrm{NP}_{\pi_{\chi_1}}(D(\pi_{\chi_1}, s))\leq H',$$ where $\mathrm{NP}$ means Newton polygon. Choose any integer $m_0$ such that $\varphi(p^{m_0})\geq\max\{\Delta+\theta^{-1},\varphi(p^m)\}$. Then, the $\pi_{\chi}$-adic Newton polygon of $D(\pi_{\chi},s)$ is invariant for $m_{\chi}\geq m_0$.
\end{lemma}

\noindent\textbf{Remark.} It is clear that $H'$ is completely determined by $H,h,h'$. It is routine to check that the middle inequality automatically holds, when $\theta>\varphi(p^m)^{-1}$.

\begin{proof}
Write that $D(T, s)=1+\sum\limits_{n\geq1}a_n(T)s^n$, where $a_n(T)=\sum\limits_ja_{n,j}T^j$. By the conditions, on every $[kh,kh+h']$, we have $$H= \mathrm{NP}_{(p^\theta, T)}(D(T, s))= \mathrm{NP}_{\pi_{\chi_1}}(D(\pi_{\chi_1}, s))= H'.$$ The equality tells us that for a vertex $(n,H(n))$ of $H$ with horizontal coordinate $n\in[kh,kh+h']$, we have
\[\nu_{(p^\theta,T)}(a_n(T))=\nu_{\pi_{\chi_1}}(a_n(\pi_{\chi_1}))=H(n).\]
Valuations of both sides are determined by finite terms. Note that
\[\nu_{(p^\theta,T)}(a_n(T))=\nu_{(p^\theta,T)}(\sum_ja_{n,j}T^j)=\min_j(\nu_{(p^\theta,T)}(a_{n,j}T^j)),\] the set
\[S=\{j:\nu_{(p^\theta,T)}(a_{n,j}T^j)=H(n)\}\] is not empty. If for every $j\in S$, $H(n)<\nu_{\pi_{\chi_1}}(a_{n,j}\pi_{\chi_1}^j),$ then by taking minimum, we have
\[H(n)<\min_{j\in S}(\nu_{\pi_{\chi_1}}(a_{n,j}\pi_{\chi_1}^j))\leq\nu_{\pi_{\chi_1}}(a_n(\pi_{\chi_1}))=H(n),\] that is a contradiction. Thus, there is some $j\in S$ such that $$H(n)=\nu_{(p^{\theta},T)}(a_{n,j}T^j)=\frac{1}{\theta}\nu_p(a_{n,j})+j\geq\nu_{\pi_{\chi_1}}(a_{n,j}\pi_{\chi_1}^j)=\varphi(p^m)\nu_p(a_{n,j})+j.$$ By the above inequality and the fact that $\theta^{-1}<\varphi(p^m)$, we have $\nu_p(a_{n,j})=0$, i.e. $a_{n,j}\in\mathbb{Z}_p^{\times}$. So any $j\in S$ satisfying the above inequality must be $H(n)$. For any $j\neq H(n)$, we have $$\nu_{\pi_{\chi_1}}(a_{n,j}\pi_{\chi_1}^{j})>H(n)=H'(n).$$ Now, consider any $\chi$ such that $m_{\chi}\geq m$. We have $$\nu_{\pi_{\chi}}(a_{n,j}\pi_{\chi}^{j})=\varphi(p^m)\nu_p(a_{n,j})+j\geq\nu_{(p^{\theta},T)}(a_{n,j}T^{j})=\frac{1}{\theta}\nu_p(a_{n,j})+j\geq H(n).$$ The two equalities can not hold at the same time, if and only if $j=H(n)$. Thus, $\nu_{\pi_{\chi}}(a_n(\pi_{\chi}))=H(n)$. As $H$ is convex, $kh,kh+h'$ are vertices, and $$\nu_{\pi_{\chi}}(a_n(\pi_{\chi}))\geq\nu_{(p^{\theta},T)}(a_n(T))\geq H(n)$$ in $n\in[kh+h',(k+1)h],$ the Newton polygon of $D(\pi_{\chi},s)$ in $[kh+h',(k+1)h]$ must be above the extension lines of the last segment of $H$ in $[kh,kh+h']$ and the first segment in $[(k+1)h,(k+1)h+h']$. Thus, the $\pi_{\chi}$-adic Newton polygon of $D(\pi_{\chi},s)$ must be coincide with $H$ in $[kh,kh+h']$ for any $k\geq0$. The proof is completed for $[kh,kh+h']$.

In $[kh+h',(k+1)h]$, $H'$ is a line, so we have the inequality similar to the condition of the lemma: $$H\leq \mathrm{NP}_{(p^\theta, T)}(D(T,s))\leq \mathrm{NP}_{\pi_{\chi}}(D(\pi_{\chi}, s))\leq H'.$$

Now we consider $n\in[kh+h',(k+1)h]$. Let $\lambda_n$ be the minimal $j$ such that $a_{n,j}\in\mathbb{Z}_p^{\times}$. If there is no such $j$, just define $\lambda_n=\infty$. For $j>\lambda_n$ and any $\chi$, we have $$\nu_{\pi_{\chi}}(a_{n,j}\pi_{\chi}^j)\geq j>\lambda_n$$. For $j<\lambda_n$, we have $\nu_p(a_{n,j})\geq1$, and $$\nu_{(p^{\theta},T)}(a_{n,j}T^j)=\frac{1}{\theta}\nu_p(a_{n,j})+j\geq H(n),$$ so for $m_{\chi}\geq m_0$,
\[\begin{split}
\nu_{\pi_{\chi}}(a_{n,j}\pi_{\chi}^j)&=\varphi(p^{m_{\chi}})\nu_p(a_{n,j})+j\geq\varphi(p^{m_0})\nu_p(a_{n,j})+H(n)-\frac{1}{\theta}\nu_p(a_{n,j})\\
&=(\varphi(p^{m_0})-\frac{1}{\theta})\nu_p(a_{n,j})+H(n)\geq\Delta+H(n)\geq H'(n).
\end{split}
\]
The set $$A=\{(n,\nu_{\pi_{\chi}}(a_n(\pi_{\chi})))|n\in[kh+h',(k+1)h],\nu_{\pi_{\chi}}(a_n(\pi_{\chi}))<H'(n)\}$$ is exactly $\{(n,\lambda_n)|n\in[kh+h',(k+1)h],\lambda_n<H'(n)\}$. But $\lambda_n$ is independent of $m_{\chi}$, so $A$ does not change as $m_{\chi}\geq m_0$ varies. By the definition of Newton polygon, the $\pi_{\chi}$-adic Newton polygon of $D(\pi_{\chi}, s)$ in $[kh+h',(k+1)h]$ is determined by $A$, thus it is also invariant.
\end{proof}

Now, we verify that the function $C^*(T,s)$ satisfies the conditions of Lemma \ref{trick}. We review the data in Section 4: for $f(x)=\sum\limits_{j=0}^{\infty}p^jf_j(x)\in\mathbb{Z}_q[x]$, $D=\mathrm{deg}(f)$, $\widetilde{d_j}=\mathrm{deg}(f_j)$, $d=\max\{\frac{\widetilde{d_j}}{p^j}|j\geq0\}$, and $l$ is the first $j$ such that $D=\widetilde{d_j}$.

Chunlei Liu and Dasheng Wei \cite{Wei} have computed the degree of $L^*$. Nonconstant terms $c_ix^i$ of $f(x)$ satisfy $(i,p)=1$, so $p\nmid\widetilde{d_j}$ for any $\widetilde{d_j}\neq0$. The non-degenerate condition in \cite{Wei} holds. So, by the main theorems of \cite{Wei}, if $\mathrm{Cond}(\chi)=p^m$, then the degree of $L^*(\chi,s)$ is
\[d_m+1=\max\limits_{0\leq j\leq m-1}\{p^{m-j-1}\widetilde{d_j}\},\] and the end point of the $q$-adic Newton polygon of $L^*(\chi,s)$ is $(d_m+1,\frac{d_m}{2})$. Observing the factor $(1 - \chi(\mathrm{Frob}_0)s)$, we know the $q$-adic Newton polygon of $L^*(\chi, s)$ has an initial segment of slope $0$. By the convexity of Newton polygons, we have
\begin{lemma}\label{upper}
The polygonal line connecting $(0,0)$, $(1,0)$ and $(d_m+1, \frac{d_m}{2})$ is an upper bound of $q$-adic Newton polygon of $L^*(\chi, s)$.
\end{lemma}

\noindent\textbf{Proof of Theorem \ref{main}}: Let $D(T,s)=C^*(T,s)$ and $m=l+1$. Choose a positive rational number $\theta$ such that $\frac{1}{\varphi(p^{l+1})}<\theta<\frac{1}{\varphi(p^l)}$. Let $H$ be the infinite convex polygon connecting $\{(n,\frac{a(p-1)n(n-1)}{2d})|n\geq0\}$. Next, set $h=dp^{l}$ and $h'=1$. Note that $H'$ is completely determined by $H$ and $h,h'$. We will check that the conditions of Lemma \ref{trick} are satisfied, with $m_0$ satisfying $\varphi(p^{m_0})\geq\max\{\frac{a(p-1)(dp^{l}-1)^2}{8dp^{l}}+\theta^{-1},\varphi(p^{l+1})\}$.

Let $\chi_1$ be a character such that $m_{\chi_1}=l+1$. By Lemma \ref{upper}, the polygonal line connecting $$(0,0),(1,0),(dp^{l},\frac{dp^{l}-1}{2})$$ is an upper bound of $q$-adic Newton polygon of $L^*(\chi_1,s)$. It begins with a segment of slope $0$, then $dp^{l}-1$ segments of slope $\frac{1}{2}$. By definition, $C^*(\chi_1,s)=\prod\limits_{j=0}^{\infty}L^*(\chi_1,q^js)$, so there is an upper bound of $q$-adic Newton polygon of $C^*(\chi_1,s)$: the polygonal line starting at $(0,0)$, then with a segment of slope $0$, then $dp^{l}-1$ segments of slope $\frac{1}{2}$, then a segment of slope $1$, then $dp^{l}-1$ segments of slope $1+\frac{1}{2}$,..., then a segment of slope $j$, then $dp^{l}-1$ segments of slope $j+\frac{1}{2}$, and so on. It is routine to check that this is the polygonal line connecting $$\{(ndp^{l},\frac{n(ndp^{l}-1)}{2}),(ndp^{l}+1,\frac{n(ndp^{l}+1)}{2})|n\geq0\}.$$ The fact $\nu_q(\pi_{\chi_1})=\frac{1}{a\varphi(p^{l+1})}$ tells us that an upper bound of the $\pi_{\chi_1}$-adic Newton polygon of $C^*(\chi_1,s)$ can be obtained by multiplying the vertical coordinates of the above vertices by $a\varphi(p^{l+1})$. This upper bound coincides with $H'$.

By Theorem \ref{lower2}, $H$ is a lower bound of the $R$(equivalent to $(p^{\theta},T)$)-adic Newton polygon of $C^*(T,s)$, where $\theta<\varphi(p^l)^{-1}$ if $l>0$. By routine calculation, we find that $\Delta=\frac{a(p-1)(dp^{l}-1)^2}{8dp^{l}}$. Let $m_0$ be any integer satisfying $\varphi(p^{m_0})\geq\max\{\frac{a(p-1)(dp^{l}-1)^2}{8dp^{l}}+\varphi(p^l),\varphi(p^{l+1})\}$. Now, Lemma \ref{trick} implies that for $m_{\chi}\geq m_0$, the $\pi_{\chi}$-adic Newton polygon of $C^*(\chi,s)$ is invariant.

By Lemma \ref{conv}, the corresponding $\mathbb{Z}_p$-tower $C_{\infty}$ is slope stable. The proof of the main theorem (Theorem \ref{main}) is completed, and the slope stable property holds for $m\geq m_0$, where $m_0$ satisfies $\varphi(p^{m_0})\geq\max\{\frac{a(p-1)(dp^{l}-1)^2}{8dp^{l}}+\varphi(p^l),\varphi(p^{l+1})\}$.\qed

\bibliographystyle{elsarticle-num}
\section*{\Large\refname}
 \bibliography{bibfile}
\end{document}